\documentclass{article}
\usepackage[english]{babel}

\usepackage[letterpaper,top=2cm,bottom=2cm,left=2.5cm,right=2.5cm,marginparwidth=1.75cm]{geometry}

\usepackage{calrsfs}
\usepackage{amssymb}
\usepackage{amsmath}
\usepackage{amsthm}
\usepackage{subcaption}
\usepackage{multirow}
\usepackage{graphicx}
\usepackage{amsfonts}
\usepackage{mathtools}
\usepackage{mathrsfs}
\usepackage{placeins}

\newtheorem{remark}{Remark}

\usepackage[mathscr]{euscript}
\usepackage[dvipsnames]{xcolor}
\usepackage[colorlinks=true, allcolors=OliveGreen]{hyperref}
\usepackage{filecontents}
\usepackage{tikz}
\usepackage{pgfplots}
\usetikzlibrary{external}
\usepackage{authblk}
\usepackage{todonotes}
\usepackage{array}
\newcolumntype{C}{>{\centering\arraybackslash}p{2.5cm}}

\newcommand{\averagel}{\{\!\!\{}
\newcommand{\averager}{\}\!\!\}}

\newcommand{\jumpl}{[\![}
\newcommand{\jumpr}{]\!]}

\newcommand{\partition}{\mathcal{T}_h}
\newcommand{\facesinternal}{\mathcal{F}^\mathrm{I}_h}
\newcommand{\faces}{\mathcal{F}_h}

\newcommand{\facesboundary}{\mathcal{F}^\mathrm{B}_h}

\newcommand{\Wh}{W_{h}^\mathrm{DG}}
\DeclareMathAlphabet{\mathcalligra}{T1}{calligra}{m}{n}

\title{Numerical modelling of protein misfolding in neurodegenerative diseases: a computational study\footnote{\textbf{Funding}: PFA has been partially supported by ICSC--Centro Nazionale di Ricerca in High Performance Computing, Big Data, and Quantum Computing funded by European Union--NextGenerationEU. PFA has been partially funded by MUR for the PRIN 2020 research grant n. 20204LN5N5. PFA and MC are members of INdAM-GNCS. The present research has been carried out within the activities of Dipartimento di Eccellenza 2023-2027.}}

\author[1]{Paola F. Antonietti\footnote{paola.antonietti@polimi.it}}

\author[1]{Mattia Corti\footnote{mattia.corti@polimi.it}}

\affil[1]{MOX-Dipartimento di Matematica, Politecnico di Milano, Piazza Leonardo da Vinci 32, Milan, 20133, Italy}

\begin{document}
\maketitle

\begin{abstract}
The spreading of misfolded proteins is a known hallmark in some neurodegenerative diseases, known as proteinopathies. A significant example is the tau protein, associated with many pathologies, such as Alzheimer's. In this work, we discuss and compare two different models for the mathematical modelling of protein misfolding, namely the heterodimer model and the Fisher-Kolmogorov model, as well as their numerical discretizations. We introduce a discontinuous Galerkin method on polygonal and polyhedral grids for space discretization to accurately simulate the wavefronts typically observed in the prionic spreading. Starting from the semidiscrete formulations, we use a Crank-Nicolson scheme to advance in time. Finally, we simulate the spreading of the misfolded tau protein in a two-dimensional brain slice in the sagittal plane with a polygonal agglomerated grid. The simulation is performed using both the presented models, and we compare the results and the differences deriving from the modelling choices.
\end{abstract} 

\section{Introduction}
\label{sec:introduction}
Neurodegenerative diseases affect millions of people worldwide, having a major impact on global health. Many pathologies, known as proteinopathies (e.g. Alzheimer's and Parkinson's diseases), have a progression which is closely linked to the aggregation of specific toxic proteins, commonly referred to as prions \cite{juckerSelfpropagationPathogenicProtein2013}. In those pathologies, an alteration in the protein's structure results in the corruption of its characteristic and functional properties, which are typically connected to resistance to the clearance mechanisms and the tendency to develop agglomerates, causing neuronal death \cite{wilsonHallmarksNeurodegenerativeDiseases2023}. It is well accepted in literature the existence of a direct correlation between the chemical microscopic reactions of prion-like proteins and the pathological advancement of dementia. However, these phenomena span over decades with a substantial delay between the protein accumulation and the clinical symptoms insurgence, making the clinical observations very complex.
\par
An important subclass of proteinopathies is the so-called tauopathies, which are associated with the misfolding and both neuronal and glial accumulation of misfolded tau protein. This microtubule protein is of primary importance in the stabilization of the neuronal cytoskeleton \cite{jouanne_tau_2017}. Pathologic tauopathies may present with diverse clinical phenotypes. These ones differentiate based on the primary site of involvement, spread and distribution. The principal pathology associated with the misfolding of tau protein is Alzheimer's disease \cite{scheltens_alzheimers_2021}. However, due to the existence of many different pathologies linked to this prionic protein, there is intense research on the construction of significant biomarkers both from the analysis of biofluid (cerebrospinal fluid and blood) and the study of different types of medical images \cite{vanoostveenImagingTechniquesAlzheimer2021}. Moreover, constructing sophisticated mathematical models for the evolution of the physical and chemical processes associated with tau protein could provide valuable insights into the disease progression \cite{carbonell_mathematical_2018}.
\par
In this work, we compare two mathematical models of decreasing mathematical complexity to describe the phenomenon: the heterodimer model and the Fisher-Kolmogorov (FK) model. The heterodimer model \cite{matthaus_diffusion_2006} aims to describe the microscopic mechanisms of prion chemical reactions while capturing the broader macroscopic propagation phenomenon. This description is possible due to the separate modelling of healthy and misfolded protein configurations. The formulation includes terms representing distinct molecular kinetic processes (production, destruction, and conversion) explicitly, making possible the development of scenario analyses \cite{fornari_prion-like_2019}. The FK model \cite{fisher-1937, kolmogorov-1937} is a nonlinear partial differential equation commonly used in the mathematical modelling of biological systems. In particular, in this context, it provides a simplified description that can be obtained from the first one under the assumption of constant healthy protein concentration \cite{weickenmeierPhysicsbasedModelExplains2019}. Indeed, the latter model is more useful whenever we do not need to analyse the dynamics of healthy protein concentration because it is almost constant. The variable of that model is a relative concentration of misfolded proteins over a maximum value that is computable from the theoretical results of the continuum heterodimer model \cite{weickenmeierPhysicsbasedModelExplains2019}.
\par 
Concerning the numerical approximation of these two models for the description of neurodegenerative disease progression, in literature, we find many different methods, such as finite element methods \cite{fornari_prion-like_2019} or reduced order network diffusion models \cite{corti_uncertainty_2023,thompsonProteinproteinInteractionsNeurodegenerative2020a}. In this work, we discuss some recently proposed numerical schemes that have shown to possess all the distinguishing features to make them suited for this problem. We consider discontinuous Galerkin formulation on polygonal/polyhedral grids (PolyDG) for the space discretization \cite{cangianiHpVersionDiscontinuousGalerkin2017}, coupled with a Crank-Nicolson time discretization scheme, with semi-implicit treatment of the nonlinear term \cite{corti_discontinuous_2023,antonietti_heterodimer_2023}. The latter approach has many advantages compared to the others. First, it allows the use of mesh composed of elements of arbitrary shape, which is particularly useful in brain geometries, where agglomeration techniques allow a detailed description of the complex domain boundaries with a relatively small number of mesh elements. Moreover, it supports high-order approximations, which are optimal for reproducing traveling-wave solutions. Finally, the use of PolyDG methods can describe correctly the heterogeneous physical properties of the domain, which is naturally subdivided into white and grey matter.
\par
This paper aims to compare the application of the two mathematical models on the same mesh and with comparable physical parameters for the first time. Moreover, we discuss the differences caused by the modelling choices and the associated computational costs. The analysis confirms that we cannot substitute in all cases the heterodimer description with a simplified FK one. Finally, for the first time, a detailed 2D simulation of the heterodimer model is made using also the heterogeneous properties of the brain matter.
\par
\bigskip
The paper is organized as follows. Section~\ref{sec:models} presents both the heterodimer and the FK models in their continuous strong formulation, and then Section~\ref{sec:polydg} is dedicated to presenting the PolyDG semi-discrete formulation. In Section~\ref{sec:fullydiscrete} we present the fully-discrete formulation of the problems. Finally, in Section~\ref{sec:results}, we simulate the spreading of tau protein within a two-dimensional brain section by comparing the use of the two methods. Finally, in Section~\ref{sec:conclusion}, we draw some conclusions and discuss future developments.

\section{Mathematical models: heterodimer and Fisher-Kolmogorov equations}
\label{sec:models}
In this section, we briefly introduce the mathematical models to describe the production, misfolding, clearance and diffusion processes of prionic proteins. For a final time $T>0$, the problem is dependent on time $t\in(0,T]$ and space $\boldsymbol{x}\in\Omega\subset\mathbb{R}^d$ ($d=2,3$). Here $\Omega$ is an open, bounded domain of $\mathbb{R}^d$. We denote by $(\cdot,\cdot)_\Omega$ the scalar product in $L^2(\Omega)$, and by $||\cdot||_\Omega$ its induced norm.
\subsection{The heterodimer model}
\label{subsec:2}
The first model we introduce is the heterodimer model, whose mathematical formulation employed to address the kinetics of prion pathogenesis incorporates prion conversion dynamics both considering monomeric and polymeric seeding hypotheses \cite{weickenmeierPhysicsbasedModelExplains2019}. The equations govern the dynamics of the overall quantity of healthy $p= p(\boldsymbol{x},t)$ and misfolded $q = q(\boldsymbol{x},t)$ proteins. Indeed, the system of equations reads:
\begin{equation}
\label{eq:hm_strong}
\begin{cases}
    \dfrac{\partial p}{\partial t}  =\nabla \cdot (\textbf{D} \nabla p) - k_{1} \, p - k_{12} \, p\, q + k_0 &  \mathrm{in} \: \Omega \times (0,T], \\[6pt]
    \dfrac{\partial q}{\partial t}  = \nabla \cdot (\textbf{D} \nabla q) - \tilde{k}_1\, q + k_{12}\, q \, p & \mathrm{in} \: \Omega \times (0,T], \\[6pt]
    (\textbf{D}\nabla p )\cdot \boldsymbol{n} = 0 \mathrm{,} \quad (\textbf{D} \nabla q )\cdot \boldsymbol{n}  = 0  & \mathrm{on} \; \partial\Omega\times (0,T], \\[6pt]
    p(\boldsymbol{x},0) = p_{0}(\boldsymbol{x}) \mathrm{,} \quad q(\boldsymbol{x},0)  = q_{0}(\boldsymbol{x}) & \mathrm{in} \: \Omega.
\end{cases}
\end{equation}
In Equation \ref{eq:hm_strong}, we assume that the brain generates only healthy protein at a generation rate governed by the parameter $k_0=k_0(\boldsymbol{x})>0$. Healthy proteins are affected by a biological destruction process at a rate $k_1=k_1(\boldsymbol{x})>0$, whereas they are converted into prions at a rate $k_{12}=k_{12}(\boldsymbol{x})>0$. In contrast, misfolded proteins undergo a clearance mechanism characterized by a rate $\tilde{k}_1=\tilde{k}_1(\boldsymbol{x})>0$. Specific relations between the parameters need to be respected to guarantee a physically consistent description of the phenomenon \cite{thompsonProteinproteinInteractionsNeurodegenerative2020a}. The spreading of proteins through the parenchyma can be mathematically characterized by employing the diffusion tensor defined as follows \cite{weickenmeierPhysicsbasedModelExplains2019}:
\begin{equation}
    \label{eq:DiffusionTensor}
    \mathbf{D} = d_\mathrm{ext}\mathbf{I} + d_\mathrm{axn}\boldsymbol{\bar{a}} \otimes \boldsymbol{\bar{a}},
\end{equation}
where the first term models the extracellular diffusion of magnitude $d_\mathrm{ext}=d_\mathrm{ext}(\boldsymbol{x})$, while the second one models the anisotropic diffusion of magnitude $d_\mathrm{axn}=d_\mathrm{axn}(\boldsymbol{x})$ along the axonal directions, denoted by the vector $\boldsymbol{\bar{a}}=\boldsymbol{\bar{a}}(\boldsymbol{x})$.
\par
Concerning the boundary conditions, assuming that the parenchyma cannot exchange the proteins with the surrounding cerebrospinal fluid, we impose homogeneous Neumann boundary conditions.
\subsection{The Fisher-Kolmogorov model}
A simplified formulation of problem~\eqref{eq:hm_strong} can be derived under the assumption that the concentration of healthy proteins is much larger than the one of misfolded proteins and that its variations can be considered negligible \cite{weickenmeierPhysicsbasedModelExplains2019, fornari_prion-like_2019}. In this case, after defining a relative misfolded protein concentration as:
\begin{equation}
    c(\boldsymbol{x},t) = \dfrac{q(\boldsymbol{x},t)}{q_{\mathrm{max}}} = \dfrac{k_{12}^2(\boldsymbol{x}) k_0(\boldsymbol{x})}{k_1^2(\boldsymbol{x})}\left(k_{12}(\boldsymbol{x})\dfrac{k_0(\boldsymbol{x})}{k_1(\boldsymbol{x})}-\tilde{k}_1(\boldsymbol{x})\right)^{-1} q(\boldsymbol{x},t),
\end{equation}
we can derive the FK model, starting from the second equation of problem~\eqref{eq:hm_strong} and using a Taylor expansion of the concentration $c$ as explained in \cite{weickenmeierPhysicsbasedModelExplains2019}. The variable $c=c(\boldsymbol{x},t)$ is a relative concentration assuming values in the interval $[0,1]$, where $0$ means the absence of misfolded proteins and $1$ is the high prevalence of them. Finally, the governing equation reduces to the FK model:
\begin{equation}
 \begin{cases}
     \dfrac{\partial c}{\partial t} =\nabla \cdot(\mathbf{D} \nabla\, c) + \alpha\,c(1-c),
    & \mathrm{in}\,\Omega\times(0,T],
    \\[8pt]
    (\mathbf{D}\nabla c) \cdot \boldsymbol{n} = 0, & 
    \mathrm{on}\;\partial \Omega \times(0,T],
    \\[8pt]
    c(\boldsymbol{x},0)=c_0(\boldsymbol{x}), & \mathrm{in}\;\Omega.
    \\[8pt]
\end{cases}
\label{eq:fk_strong}
\end{equation}
where $\alpha = k_{12}\frac{k_0}{k_1}-\tilde{k}_1$ is the conversion rate, describing with a single term both misfolding, clearance and protein production processes. Concerning the diffusion tensor $\mathbf{D}$ it is defined in Equation~\eqref{eq:DiffusionTensor}. The homogeneous Neumann boundary condition extends naturally to the concentration $c$ we introduce.
\begin{remark}
    Due to the simplifications we need to derive the FK model, we remark that it would be usable only in some conditions, which are not always respected in the reality of the applications. However, the completeness of the heterodimer model is associated with a higher computational cost, due to the need to solve the problem for two variables. 
\end{remark}
\section{Polytopal discontinuous Galerkin semi-discrete formulation}
\label{sec:polydg}
In this section, after defining some preliminary concepts, we introduce the semi-discrete formulation in space of the two proposed models of~\eqref{eq:hm_dgformulation} and~\eqref{eq:fk_strong} using the PolyDG method. We introduce a polytopic mesh partition $\partition$ of the domain $\Omega$ made of disjoint polygonal/polyhedral elements $K$. For each element $K\in \partition$, we denote by $|K|$ the measure of the element, by $h_K$ its diameter, and set $h=\max_{K\in\partition} h_K$. We define the interface as the intersection of the $(d-1)-$dimensional facets of two neighbouring elements and the set of interfaces as $\faces$. We distinguish two cases:
\begin{itemize}
    \item case $d=2$, in which the interfaces are always line segments;
    \item case $d=3$, in which any interface consists of a generic polygon, we can decompose each interface into (planar) triangles that compose $\faces$.
\end{itemize}
It is now useful to decompose $\faces$ into the union of interior faces ($\facesinternal$) and exterior faces ($\facesboundary$) lying on the boundary of the domain $\partial\Omega$, i.e. $\faces = \facesinternal \cup \facesboundary$. We underline that in this article the boundary faces are all associated with Neumann boundary conditions and the formulation is constructed accordingly. Concerning the assumptions on the domain partition we refer to the properties in \cite{corti_discontinuous_2023}. 
\par
Let us define $\mathbb{P}_{\ell}(K)$ as the space of polynomials of total degree $\ell \geq 1$ over a mesh element $K$. Then we can introduce the following discontinuous finite element space:
\begin{equation*}
    \Wh = \{w\in L^2(\Omega):\quad w|_K\in\mathbb{P}_{\ell}(K)\quad\forall K\in\partition\}
\end{equation*}
\par
Finally, we need to introduce the trace operators \cite{arnoldUnifiedAnalysisDiscontinuous2001}. Let $F\in\facesinternal$ be a face shared by the elements $K^\pm$ and let $\boldsymbol{n}^\pm$ be the unit normal vector on face $F$ pointing exterior to $K^\pm$, respectively. Then, for sufficiently regular scalar-valued functions $v$ and vector-valued functions $\boldsymbol{q}$, respectively, we define:
\begin{itemize}
    \item average operator $\averagel{\cdot}\averager$ on $F\in\facesinternal$: $\averagel{v}\averager = \dfrac{1}{2} (v^+ + v^-), \quad \averagel{\boldsymbol{q}}\averager = \dfrac{1}{2} (\boldsymbol{q}^+ + \boldsymbol{q}^-)$;
    \item jump operator $\jumpl{\cdot}\jumpr$ on $F\in\facesinternal$: $\jumpl{v}\jumpr = v^+\boldsymbol{n}^+ + v^-\boldsymbol{n}^-, \quad \jumpl{\boldsymbol{q}}\jumpr = \boldsymbol{q}^+\cdot\boldsymbol{n}^+ + \boldsymbol{q}^-\cdot\boldsymbol{n}$.
\end{itemize}
In these relations, the superscripts $\pm$ on the functions denote the traces of the functions on $F$ within the interior to $K^\pm$. 
\subsection{PolyDG semi-discrete approximation of the heterodimer model}
To construct the semi-discrete formulation, we define the penalization function
\begin{equation}
  \eta:\facesinternal\rightarrow\mathbb{R}_+ \mathrm{\;such\,that\;}\eta = \eta_0 
    \max\left\{\{d^K\}_\mathrm{H},\{k^K\}_\mathrm{H}\right\}\dfrac{p^2}{\{h\}_\mathrm{H}},
    \label{eq:hm_penalty}
\end{equation}
where $\eta_0$ is a constant parameter that should be chosen sufficiently large to ensure the stability of the discrete formulation, $d^K = \|\sqrt{\mathbf{D}|_K}\|^2$, and $k^K = \|\,(1 + k_{12}|_K)(k_1|_K + \tilde{k}_1|_K)\|$. In Equation \eqref{eq:hm_penalty}, we are considering the harmonic average operator defined as $\{v\}_\mathrm{H} = \frac{2 v^+ v^-}{v^+ + v^-}$. We define the bilinear form $\mathcal{A}_h:\Wh\times \Wh\rightarrow \mathbb{R}$ as:
\begin{equation}
\begin{split}
        \mathcal{A}_h(c,w) = \int_{\Omega} \left(\mathbf{D} \nabla_h v\right)\cdot\nabla_h w - & \sum_{F\in\facesinternal}\int_{F}\left(\averagel\mathbf{D} \nabla_h v\averager \cdot \jumpl w \jumpr +  \jumpl v\jumpr \cdot \averagel\mathbf{D} \nabla_h w \averager\right) + \\
        + & \sum_{F\in\facesinternal}\int_{F}\eta \jumpl v\jumpr \cdot \jumpl w\jumpr \mathrm{d}\sigma \qquad \forall v,w \in\Wh,
\end{split}
\label{eq:bilinearform}
\end{equation}
where $\nabla_h$ is the elementwise gradient. Given suitable discrete approximations $p_{0h},q_{0h}\in \Wh$ of the initial conditions of Equation~\eqref{eq:hm_strong}, the semi-discrete PolyDG formulation reads:
\par
\bigskip
For each $t>0$, find $(p_h,q_h) =(p_h(t),q_h(t))\in \Wh\times\Wh$ such that:
\begin{equation}
\begin{cases}
  \left(\dfrac{\partial p_h}{\partial t}, v_h\right)_{\Omega} + \mathcal{A}_h(p_h,v_h) + (k_1\,p_h+k_{12}\,p_h^2,v_h)_\Omega = (k_0,v_h)_\Omega & \forall v_h\in \Wh, \\
 \left(\dfrac{\partial q_h}{\partial t}, w_h\right)_{\Omega} + \mathcal{A}_h(q_h,w_h) + (\tilde{k}_1\,q_h-k_{12}\,q_h^2,w_h)_\Omega = 0 & \forall w_h\in \Wh, \\
 p_h(\boldsymbol{x},0) = p_{0h}, \qquad q_h(\boldsymbol{x},0) = q_{0h}. 
\end{cases}
\label{eq:hm_dgformulation}
\end{equation}
For details regarding the derivation of the semi-discrete formulation and the stability and \textit{a-priori} error analysis, we refer to \cite{antonietti_heterodimer_2023}.
\subsection{PolyDG semi-discrete approximation of the Fisher-Kolmogorov model}
To construct the semi-discrete formulation, we define the penalization function
\begin{equation}
     \eta:\facesinternal\rightarrow\mathbb{R}_+ \mathrm{\;such\,that\;}\eta = \eta_0 \max\left\{\{d^K\}_\mathrm{H},\{\alpha^K\}_\mathrm{H}\right\}\dfrac{p^2}{\{h\}_\mathrm{H}},
    \label{eq:fk_penalty}
\end{equation}
where $\eta_0$ is a constant parameter that should be chosen sufficiently large to ensure the stability of the discrete formulation, $d^K$ is defined as before, and $\alpha^K = \|\alpha|_K\|$. 
\par
Using the definition of the penalty parameter in Equation~\eqref{eq:fk_penalty} in the definition of the bilinear form in Equation~\eqref{eq:bilinearform}, the semi-discrete PolyDG formulation reads:
\par
\bigskip
For each $t>0$, find $c_h=c_h(t)\in \Wh$ such that $\forall t>0$:
\begin{equation}
\label{eq:fk_dgformulation}
    \begin{cases}
     \left(\dfrac{\partial c_h}{\partial t},w_h\right)_\Omega + \mathcal{A}_h(c_h,w_h) - (\alpha c_h,w_h)_\Omega + (\alpha c_h^2,w_h)_\Omega = 0
     \quad \forall w_h\in \Wh, \\
    c_h(\boldsymbol{x},0) = c_{0h}.
    \end{cases}
\end{equation}
where $c_{0h}\in\Wh$ is a suitable approximation of $c_0$. For details regarding the derivation of the semi-discrete formulation and the stability and \textit{a-priori} error analysis, we refer to \cite{corti_discontinuous_2023}.
\begin{remark}
In both formulations, the dependency of the penalty coefficient expression on the reaction coefficients (see Equation~\eqref{eq:hm_penalty} and~\eqref{eq:fk_penalty}) is justified by the \textit{a-priori} error estimates that provide relations between the coercivity constant of the bilinear form in Equation~\eqref{eq:bilinearform} and the reaction parameters to be respected. More details can be found in \cite{antonietti_heterodimer_2023, corti_discontinuous_2023}.
\end{remark}
\begin{remark}
The construction of the PolyDG discretization can be extended to more general boundary conditions and forcing terms, with minimal adjustments, for both formulations. More details can be found in \cite{antonietti_heterodimer_2023, corti_discontinuous_2023}, where the analysis is performed in a more general setting. The choice of using a simplified setting in this work is justified by the application we focus on.
\end{remark}
\section{Fully-discrete formulation}
\label{sec:fullydiscrete}
First of all, we need to construct the vectors and matrices we need for the algebraic formulation. We consider $\{\phi_j\}_{j=0}^{N_h}$ a set of basis function for $\Wh$, being $N_h$ its dimension. We can write the discrete solutions expansion in terms of the chosen basis:
\begin{equation*}
    p_h(\boldsymbol{x},t) = \sum_{j=0}^{N_h} \mathbf{P}_j(t) \phi_j(\boldsymbol{x}) \quad
    q_h(\boldsymbol{x},t) = \sum_{j=0}^{N_h} \mathbf{Q}_j(t) \phi_j(\boldsymbol{x}) \quad
    c_h(\boldsymbol{x},t) = \sum_{j=0}^{N_h} \mathbf{C}_j(t) \phi_j(\boldsymbol{x}).
\end{equation*} 
We denote by $\mathbf{P}(t)$, $\mathbf{Q}(t)$, $\mathbf{C}(t)\in\mathbb{R}^{N_h}$, the corresponding vectors of the expansion coefficients. Moreover, we define the following matrices $\mathrm{M},\mathrm{A},\mathrm{M}_\omega,\widehat{\mathrm{M}}_\omega\in\mathbb{R}^{N_h\times N_h}$ as:
\begin{align*}
     [\mathrm{M}]_{i,j}                 & = (\phi_j,\phi_i)_\Omega          & \mathrm{Mass\;matrix} \\
     [\mathrm{A}]_{i,j}                 & = \mathcal{A}_h(\phi_j,\phi_i)      & \mathrm{Stiffness\;matrix} \\
     [\mathrm{M}_\omega]_{i,j}               & = (\omega\phi_j,\phi_i)_\Omega              & \mathrm{Linear\;reaction\;matrix\;of\;a\;parameter\;}\omega \\
     [\widehat{\mathrm{M}}_\omega (\Phi)]_{i,j}               & = (\omega\Phi_h\phi_j,\phi_i)_\Omega              & \mathrm{Non-linear\;reaction\;matrix\;of\;a\;parameter\;}\omega \\
\end{align*}
with $i,j = 1,...,N_h$ and $\omega=\{\alpha,k_1,\tilde{k}_1,k_{12}\}$. Finally, we introduce the right-hand side vector $[\mathbf{F}_{k_0}]_i= (k_0,\phi_i)_\Omega$  with $i = 1,...,N_h$. 
\par
Concerning time discretization, we adopt a time-stepping method of the second order using a Crank-Nicolson scheme. We construct a partition $0 < t_1 < t_2 < ... < t_{N_T} = T$ of the time interval $[0,T]$ into $N_T$ intervals of constant time step $\Delta t = t_{n+1}-t_n$. 
\subsection{Fully-discrete formulation of heterodimer model}
First of all, applying the above definitions we recall the semi-discrete algebraic formulation of~\eqref{eq:hm_dgformulation}:
\begin{equation}
    \label{eq:HM_AlgebraicFormulation}
    \begin{cases}
    \mathrm{M}\dot{\mathbf{P}}(t) + \mathrm{A}\mathbf{P}(t) + \mathrm{M}_{k_1}\mathbf{P}(t) + \widehat{\mathrm{M}}_{k_{12}}\left(\mathbf{P}(t)\right) \mathbf{Q}(t) =  \mathbf{F}_{k_0}, &  t\in (0,T] \\[6pt] 
   \mathrm{M}\dot{\mathbf{Q}}(t) + \mathrm{A}\mathbf{Q}(t) + \mathrm{M}_{\tilde{k}_1}\mathbf{Q}(t) - \widehat{\mathrm{M}}_{k_{12}}\left(\mathbf{Q}(t)\right) \mathbf{P}(t) =  \mathbf{0}, &  t\in (0,T] \\[6pt]
   \mathbf{P}(0) = \mathbf{P}_0, \qquad \mathbf{Q}(0) =\mathbf{Q}_0,
    \end{cases}
\end{equation}
where $\mathbf{P}_0$ and $\mathbf{Q}_0$ are the vector expansions associated to the initial conditions $p_{0h}$ and $q_{0h}$, respectively. Then, we present the fully-discrete formulation of the heterodimer model. For $n= 1, ... , N_T$ we want to find $(\mathbf{P}^{n}, \mathbf{Q}^{n})$ solution of the following system of equations:
\begin{equation}
    \label{eq:HM_FullyDiscreteFormulation}
    \begin{cases}
    \left(\dfrac{\mathrm{M}}{\Delta t} + \dfrac{1}{2}(\mathrm{A}+\mathrm{M}_{k_1})\right)\mathbf{P}^{n} + \dfrac{1}{2}\widehat{\mathrm{M}}_{k_{12}}\left(\frac{3}{2}\mathbf{P}^{n-1}-\frac{1}{2}\mathbf{P}^{n-2}\right) \mathbf{Q}^n =  \mathbf{F}_p^{n}, \\  
    \left(\dfrac{\mathrm{M}}{\Delta t} + \dfrac{1}{2}(\mathrm{A}+\mathrm{M}_{\tilde{k}_1})\right)\mathbf{Q}^{n} - \dfrac{1}{2}\widehat{\mathrm{M}}_{k_{12}}\left(\frac{3}{2}\mathbf{Q}^{n-1}-\frac{1}{2}\mathbf{Q}^{n-2}\right) \mathbf{P}^n =  \mathbf{F}_q^{n}, \\[6pt]
   \mathbf{P}^0 = \mathbf{P}_0, \qquad \mathbf{Q}^0 =\mathbf{Q}_0,
    \end{cases}
\end{equation}
where: 
\begin{align*}
    \mathbf{F}_p^{n} = & \left(\dfrac{\mathrm{M}}{\Delta t} - \dfrac{1}{2}(\mathrm{A}+\mathrm{M}_{k_1})\right)\mathbf{P}^{n-1} - \dfrac{1}{2}\widehat{\mathrm{M}}_{k_{12}}\left(\frac{3}{2}\mathbf{P}^{n-1}-\frac{1}{2}\mathbf{P}^{n-2}\right) \mathbf{Q}^{n-1}+\mathbf{F}_{k_0}^{n}, \\ 
    \mathbf{F}_q^{n} = & \left(\dfrac{\mathrm{M}}{\Delta t} - \dfrac{1}{2}(\mathrm{A}+\mathrm{M}_{\tilde{k}_1})\right)\mathbf{P}^{n-1} + \dfrac{1}{2}\widehat{\mathrm{M}}_{k_{12}}\left(\frac{3}{2}\mathbf{Q}^{n-1}-\frac{1}{2}\mathbf{Q}^{n-2}\right) \mathbf{P}^{n-1}.
\end{align*}
A discussion of the performance of the Crank-Nicolson scheme versus an implicit Euler time integration method can be found in \cite{antonietti_heterodimer_2023}.
\begin{figure}[t]
    \centering
        \includegraphics[width=\textwidth]{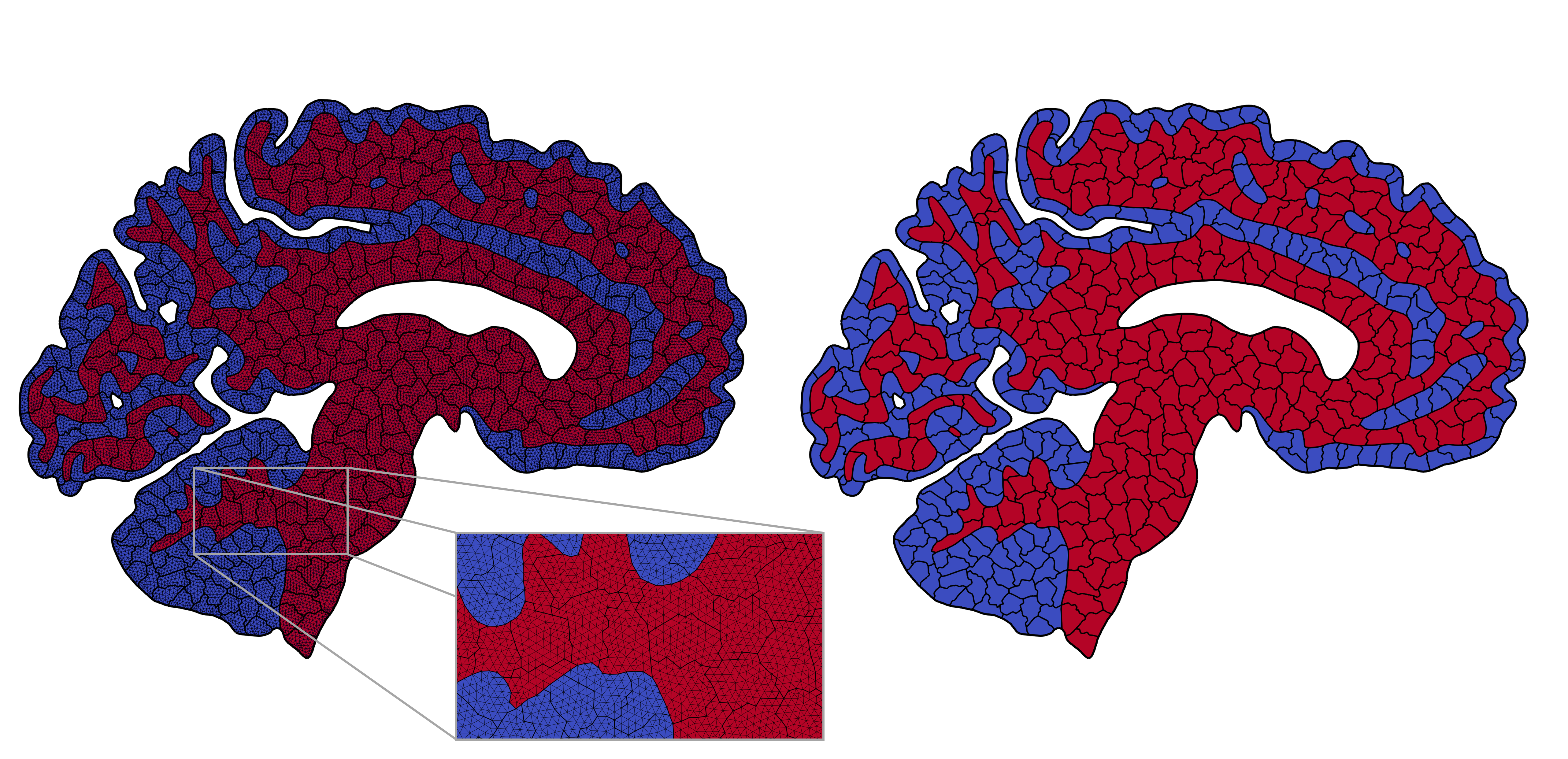}
   \caption{Agglomerated mesh with the distinction of white matter (red) and grey matter (blue). The initial fine triangular mesh is visible in the left image with a zoom of a particular region.}
   \label{fig:BrainMesh}
\end{figure}
\subsection{Fully-discrete formulation of Fisher-Kolmogorov model}
Concerning the FK model in~\eqref{eq:fk_strong}, the semi-discrete algebraic formulation of~\eqref{eq:fk_dgformulation} reads:
\begin{equation}
    \label{eq:FK_AlgebraicFormulation}
    \begin{cases}
   \mathrm{M}\dot{\mathbf{C}}(t) + \mathrm{A}\mathbf{C}(t) - \mathrm{M}_{\alpha}\mathbf{C}(t) + \widehat{\mathrm{M}}_{\alpha}\left(\mathbf{C}(t)\right) \mathbf{C}(t) =  \mathbf{0}, &  t\in (0,T] \\[6pt]
   \mathbf{C}(0) = \mathbf{C}_0,
    \end{cases}
\end{equation}
where $\mathbf{C}_0$ is the vector expansion associated to the initial conditions $c_{0h}$. Then we derive the fully-discrete formulation. For $n= 1, ... , N_T$ we want to find $\mathbf{C}^{n}$ solution of the following linear system of equations:
\begin{equation}
    \label{eq:FK_FullyDiscreteFormulation}
    \left(\dfrac{\mathrm{M}}{\Delta t} + \dfrac{1}{2}(\mathrm{A}-\mathrm{M}_{\alpha})\right)\mathbf{C}^{n} + \dfrac{1}{2}\widehat{\mathrm{M}}_{\alpha}\left(\frac{3}{2}\mathbf{C}^{n-1}-\frac{1}{2}\mathbf{C}^{n-2}\right) \mathbf{C}^n =  \mathbf{F}_c^{n}, 
\end{equation}
where: 
\begin{align*}
    \mathbf{F}_c^{n} = & \left(\dfrac{\mathrm{M}}{\Delta t} - \dfrac{1}{2}(\mathrm{A}-\mathrm{M}_\alpha)\right)\mathbf{C}^{n-1} - \dfrac{1}{2}\widehat{\mathrm{M}}_{\alpha}\left(\frac{3}{2}\mathbf{C}^{n-1}-\frac{1}{2}\mathbf{C}^{n-2}\right) \mathbf{C}^{n-1}.
\end{align*}
More details about the construction of this formulation and a comparison with the use of Picard iterations to treat the nonlinear term can be found in \cite{corti_discontinuous_2023}.

\section{Numerical results: tau protein spreading in a brain section}
\label{sec:results}

In this section, we perform some numerical simulations about the spreading of the tau protein in its misfolded form in Alzheimer's disease to compare the predictive capabilities and performances of the models under investigation. The numerical simulations have been obtained on the open source Lymph library \cite{antonietti_lymph_2024}, implementing the PolyDG method for multiphysics. For the simulations, we consider a detailed mesh of the brain obtained by segmentation of a structural MRI dataset sourced from the OASIS-3 database \cite{OASIS3}, as described in \cite{corti_structure_preserving_2023}. The segmentation is performed using Freesurfer \cite{Freesurfer}.
\par
The sagittal section of the brain is then used to construct a fine triangulated mesh composed of $43\,402$ triangles. An agglomerated polygonal mesh of $534$ elements is constructed starting from the triangular one, based on employing Parmetis \cite{Parmetis}. The versatility of the PolyDG method allows us to reduce the computational cost of the simulation using high-order approximation on a coarse mesh. Moreover, the advantage of using an agglomeration technique is that we preserve the level of detail of the underlying finer mesh in the description of the domain boundary. The agglomerated mesh with the detail of the triangular one used for the agglomeration is visible in Figure~\ref{fig:BrainMesh}. In this case, the agglomeration is performed maintaining the distinction between white and grey matters. This construction allows the specialization of the physical parameters in the two models, considering the heterogeneity of the brain structure \cite{antonietti_agglomeration_2024}. 
\par
To construct the axonal component of the diffusion tensor $\mathbf{D}$, we derive the axonal directions starting from the Diffusion-Weighted Image (DWI) and computing the principal eigenvector of the diffusion tensor in each voxel of the image. For more detail, we refer to \cite{Mardal_Mesh}. 
\par
Concerning the discretization parameters in space, we fix the polynomial order $\ell=6$ in every element of the discretization and we set the penalty parameter $\gamma_0=10$. Regarding the time discretization, we adopt a time step $\Delta t = 0.01$ and simulate up to $T=25$ years.
\begin{figure}[t]
    \centering
    \includegraphics[width=0.9\textwidth]{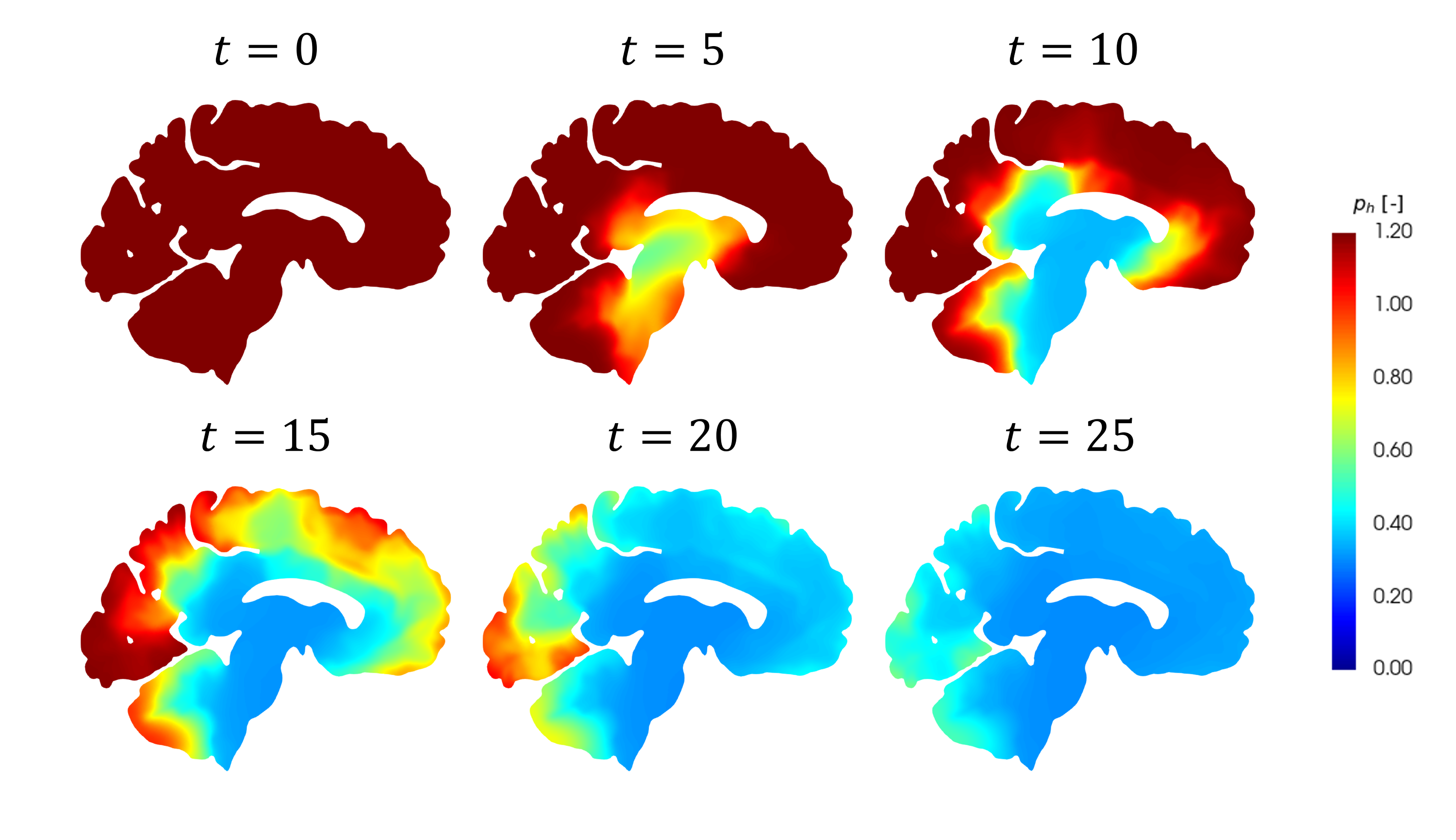}
   \caption{The initial configuration corresponds to the unstable healthy state, with a protein concentration of $p_h=1.2$ uniformly distributed all over the brain. The triggering of the system is due to the presence of misfolded proteins.}
   \label{fig:solution_brain_p}
\end{figure}

\begin{figure}[t]
    \centering
    \includegraphics[width=0.9\textwidth]{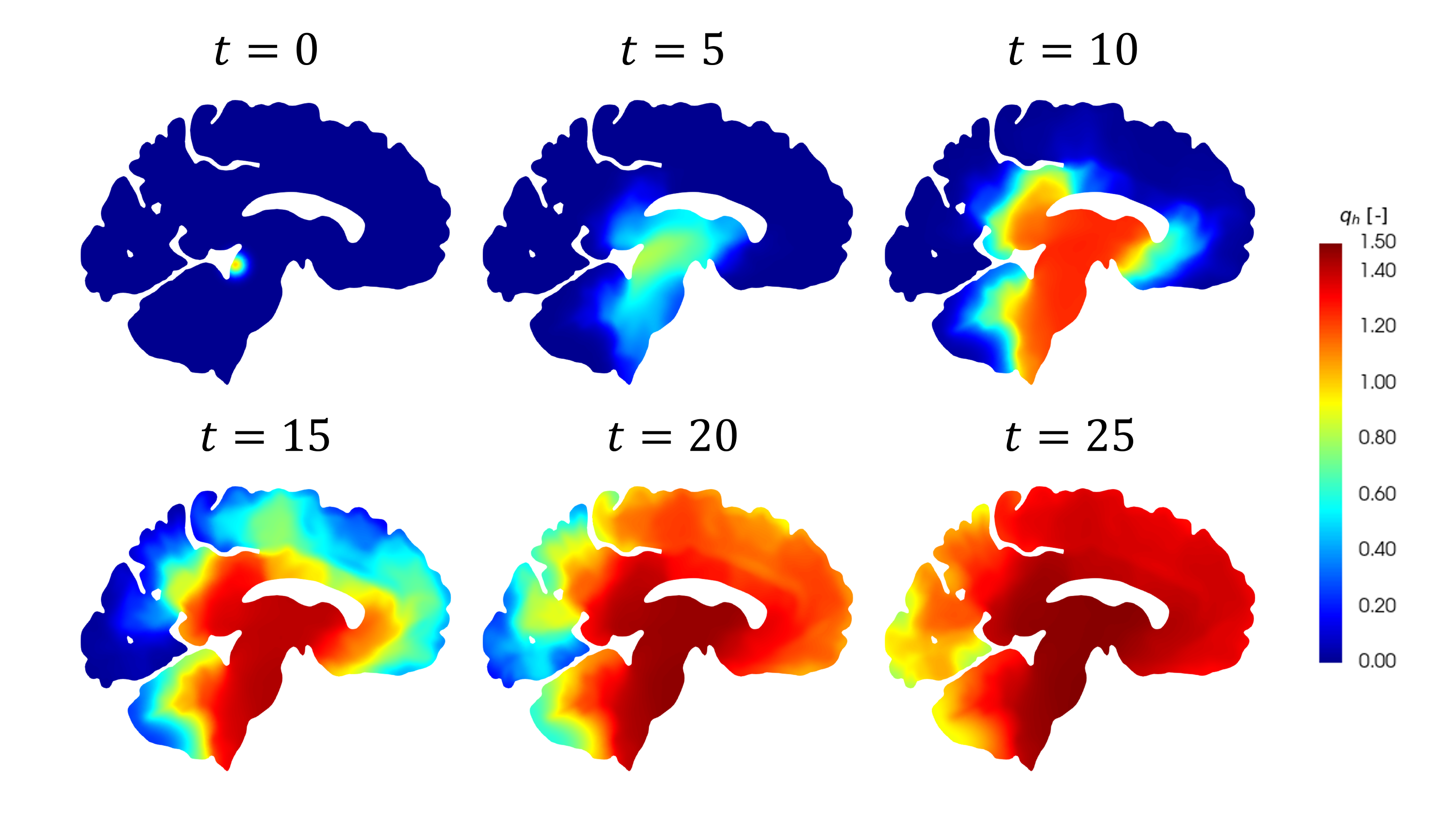}
   \caption{Misfolded protein $q_h$, initially concentrated into the entorhinal cortex, spreads through the entire brain section over 25 years.}
   \label{fig:solution_brain_q}
\end{figure}

\begin{table}[t]
	\centering
	\begin{tabular}{|c|l l|l l|}
	\hline
	\textbf{Model parameter} & \multicolumn{2}{c|}{\textbf{Value in white matter}} & \multicolumn{2}{c|}{\textbf{Value in grey matter}} \\ 
		\hline 
		 $d_\mathrm{ext}$  & $\quad8.0\times10^{-6}$ & $[\mathrm{m^2/years}]\quad$ & $\quad8.0\times10^{-6}$ & $[\mathrm{m^2/years}]\quad$ \\
		 $d_\mathrm{axn}$  & $\quad8.0\times10^{-5}$ & $[\mathrm{m^2/years}]\quad$ & $\quad0.0\times10^{-6}$ & $[\mathrm{m^2/years}]\quad$ \\ 
		 $k_0$             & $\quad6.0\times10^{-1}$ & $[\mathrm{1/years}]\quad$    & $\quad6.7\times10^{-2}$ & $[\mathrm{1/years}]\quad$ \\
		 $k_{12}$          & $\quad1.0\times10^{0}$  & $[\mathrm{1/years}]\quad$    & $\quad1.1\times10^{-1}$ & $[\mathrm{1/years}]\quad$ \\ 
		 $k_1$             & $\quad5.0\times10^{-1}$ & $[\mathrm{1/years}]\quad$    & $\quad5.6\times10^{-2}$ & $[\mathrm{1/years}]\quad$ \\
		 $\widetilde{k}_1$ & $\quad3.0\times10^{-1}$ & $[\mathrm{1/years}]\quad$    & $\quad3.3\times10^{-2}$ & $[\mathrm{1/years}]\quad$ \\
		\hline 
	\end{tabular}
	\caption{Physical parameter values used in the numerical simulation of the heterodimer model \cite{thompsonProteinproteinInteractionsNeurodegenerative2020a,antonietti_heterodimer_2023}.}
	\label{tab:hm_param}
\end{table}

\subsection{Heterodimer model}
\label{BrainApp}
This section is devoted to the simulation of the spread of tau protein into a two-dimensional section of the human brain using the heterodimer model. The diffusion modelling required appropriate parameter configurations. We report the used numerical values of the physical parameters in Table~\ref{tab:hm_param}. As done in \cite{weickenmeierPhysicsbasedModelExplains2019}, we adopt the assumption of fast misfolding and preferential diffusion in the white matter. The values in the white matter are coherent with the literature works on topic \cite{thompsonProteinproteinInteractionsNeurodegenerative2020a}, while the values in the grey one are coherently computed to obtain the same equilibrium concentration with reduced velocity of the misfolding process.
\par
By construction we expect the system to develop from the initial state towards the stable equilibrium point $(p_h,q_h)=(0.3,1.5)$. Moreover, we set an initial condition for the healthy protein concentration uniformly equal to the unstable equilibrium point $p_h=1.2$ as in \cite{antonietti_heterodimer_2023}. Finally, we trigger the system with an initial seeding of the misfolded version of tau protein concentrated in the entorhinal cortex as in \cite{weickenmeierPhysicsbasedModelExplains2019}.
\par
In Figures \ref{fig:solution_brain_p} and \ref{fig:solution_brain_q}, we report the approximated solution at different time instants ($t=0,5,10,15,20,25$ years). As we can observe, the solution develops to the pathological diffusion of the misfolded tau protein in the whole brain geometry in $25$ years. The results are coherent with what is expected from the Braak staging theory for this type of protein \cite{jouanne_tau_2017,braak_alzheimer_1991}. Moreover, it is observable how the differences in the physical parameters cause a slower increase in the concentration of misfolded proteins in the grey matter, which is both associated with reduced reactivity and diffusion.
\subsection{Fisher-Kolmogorov model}
\begin{table}[t]
	\centering
	\begin{tabular}{|c|l l|l l|}
	\hline
	\textbf{Model parameter} & \multicolumn{2}{c|}{\textbf{Value in white matter}} & \multicolumn{2}{c|}{\textbf{Value in grey matter}} \\ 
		\hline 
		 $d_\mathrm{ext}$  & $\quad8.0\times10^{-6}$ & $[\mathrm{m^2/years}]\quad$ & $\quad8.0\times10^{-6}$ & $[\mathrm{m^2/years}]\quad$ \\
		 $d_\mathrm{axn}$  & $\quad8.0\times10^{-5}$ & $[\mathrm{m^2/years}]\quad$ & $\quad0.0\times10^{-6}$ & $[\mathrm{m^2/years}]\quad$ \\  
		 $\alpha$             & $\quad9.0\times10^{-1}$ & $[\mathrm{1/years}]\quad$    & $\quad1.0\times10^{-2}$ & $[\mathrm{1/years}]\quad$ \\
		\hline 
	\end{tabular}
	\caption{Physical parameter values used in the numerical simulation of heterodimer model \cite{weickenmeierPhysicsbasedModelExplains2019, corti_discontinuous_2023}.}
	\label{tab:fk_param}
\end{table}
\begin{figure}[t]
    \centering
    \includegraphics[width=0.9\textwidth]{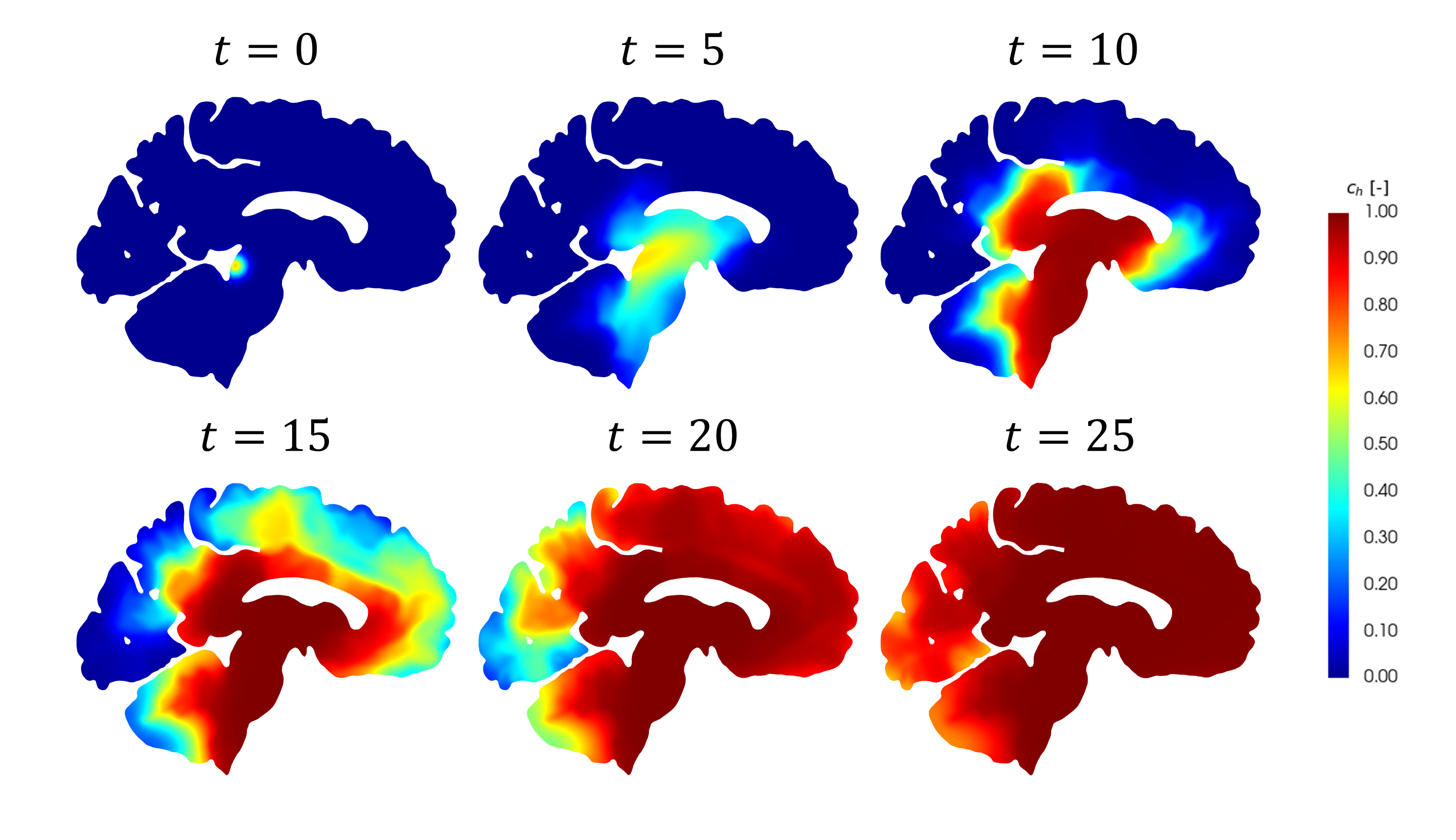}
   \caption{Relative misfolded protein $c_h$ concentration, initially concentrated into the entorhinal cortex, spreads through the entire brain section.}
   \label{fig:solution_brain_c}
\end{figure}
In this section, we adopt the simplified FK modelling for the simulation of the spreading of tau protein, comparing the results with the ones obtained with the heterodimer model, starting from comparable initial conditions. Indeed, also in this case we trigger the system with the same seeding of the misfolded version of tau protein in the heterodimer simulation rescaled over the value $q_\mathrm{max}=1.5$. Moreover, we impose the values of the parameter $\alpha$ in the two matters, computing it from the numerical values adopted in the previous simulation (see Table~\ref{tab:hm_param}) and exploiting the relation $\alpha = k_{12}\frac{k_0}{k_1}-\tilde{k}_1$. The diffusion values are the same as in the heterodimer model simulation. We report the numerical values in Table~\ref{tab:fk_param}. 
\par
By construction we expect the system to develop from the initial state towards the stable equilibrium point $c_h = 1$. Using the FK model, the results are coherent with the Braak stages in the spreading of the tau protein \cite{jouanne_tau_2017,braak_alzheimer_1991}. However, the use of a simplified model causes a slight change in the final result, with respect to the heterodimer simulation. Indeed, by making a comparison between Figure~\ref{fig:solution_brain_q} and Figure~\ref{fig:solution_brain_c}, we can observe that in the FK simulation, the velocity of the spreading is faster, due to the assumptions made on the concentration of healthy proteins $p$ in the model derivation. This is because the heterodimer simulation does not respect the assumption $p\gg q$ and this creates significant dynamics in the healthy protein concentration, which variations slow down the diffusion of the pathology.
\begin{figure}[t]
     \begin{subfigure}[b]{0.48\textwidth}
         \centering
         \includegraphics[width=0.9\textwidth]{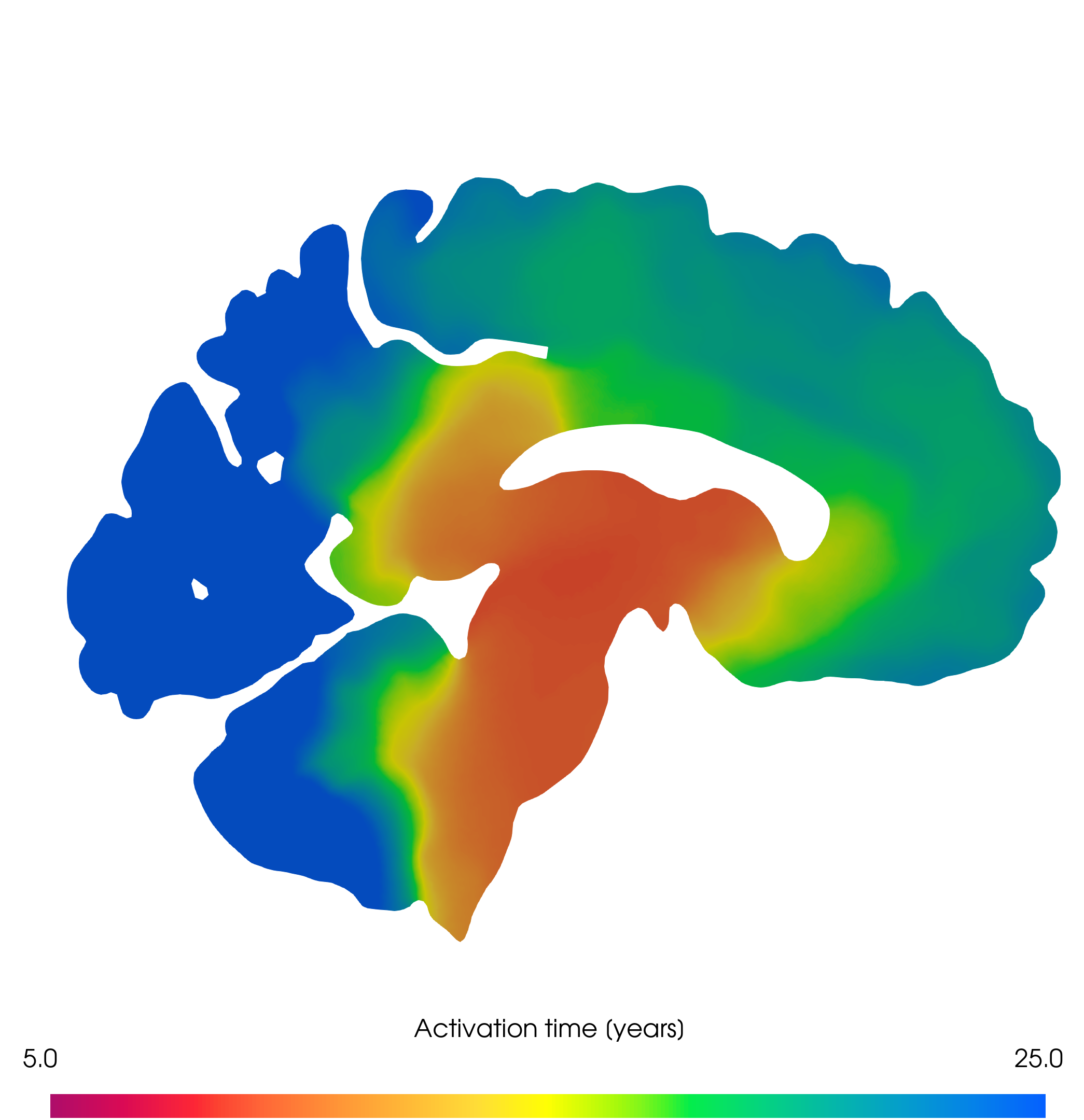}
         \caption{Activation time using heterodimer model.}
         \label{fig:HM_activationtime}
     \end{subfigure}
          \begin{subfigure}[b]{0.48\textwidth}
         \centering
         \includegraphics[width=0.9\textwidth]{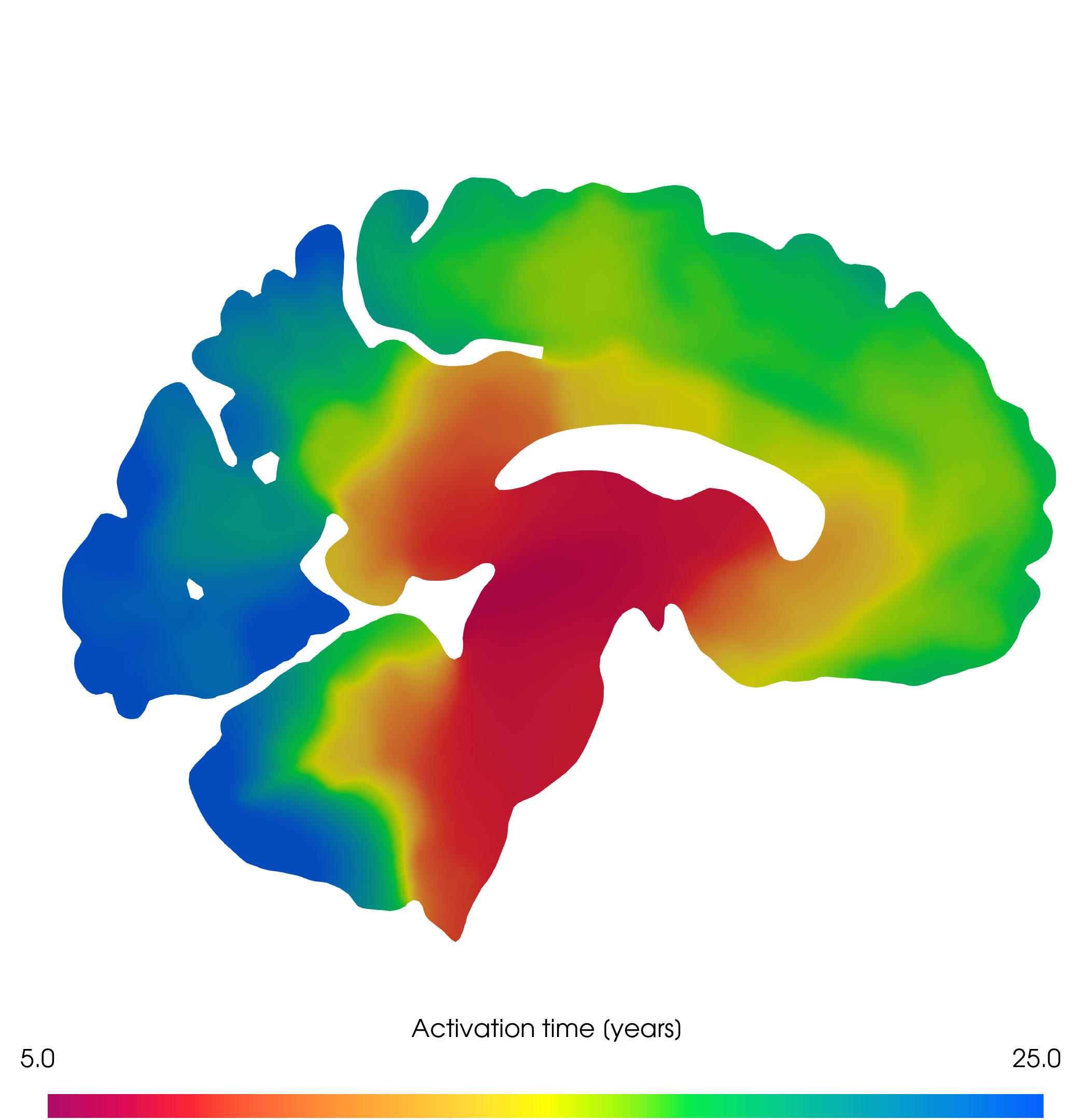}
         \caption{Activation time using FK model.}
         \label{fig:FK_activationtime}
     \end{subfigure}
     \caption{Activation time computed starting from the simulations of the two different models: heterodimer (a) and FK (b).}
    \label{fig:activationtime}
\end{figure}
\par
To compare with more detail the impact of the two different models, we compute the activation time \cite{weickenmeierPhysicsbasedModelExplains2019, corti_discontinuous_2023} of the pathology as:
\begin{equation}
\label{eq:acttime}
    \widehat{t}(\boldsymbol{x},t) = \chi_{\{v_h(\boldsymbol{x},t)>v_\mathrm{crit}\}} (\boldsymbol{x},t) \qquad \boldsymbol{x}\in \Omega \quad t\in[0,T],
\end{equation}
where $\chi$ is the indicator function and $v_\mathrm{crit}$ is the critical value of the pathological protein concentration we fix to be equal to $1.20$ for the heterodimer model ($v_h=q_h$) and $0.80$ for FK one ($v_h=c_h$), to obtain comparable results. The indicator gives us a measure of the time after which the neurons in a specific region will be affected by pathological communication. We report the activation time computed starting from the simulation of the heterodimer model in Figure~\ref{fig:HM_activationtime} and the one from the FK one in Figure~\ref{fig:FK_activationtime}. Coherently to what we discuss from the solution analysis the activation times are generally faster in the FK simulation. Moreover, in both cases, we can notice that the white matter activates before the grey one, due to the assumptions made on the physical parameters value.
\par
\begin{figure}[t]
    \centering
    \includegraphics[width=0.9\textwidth]{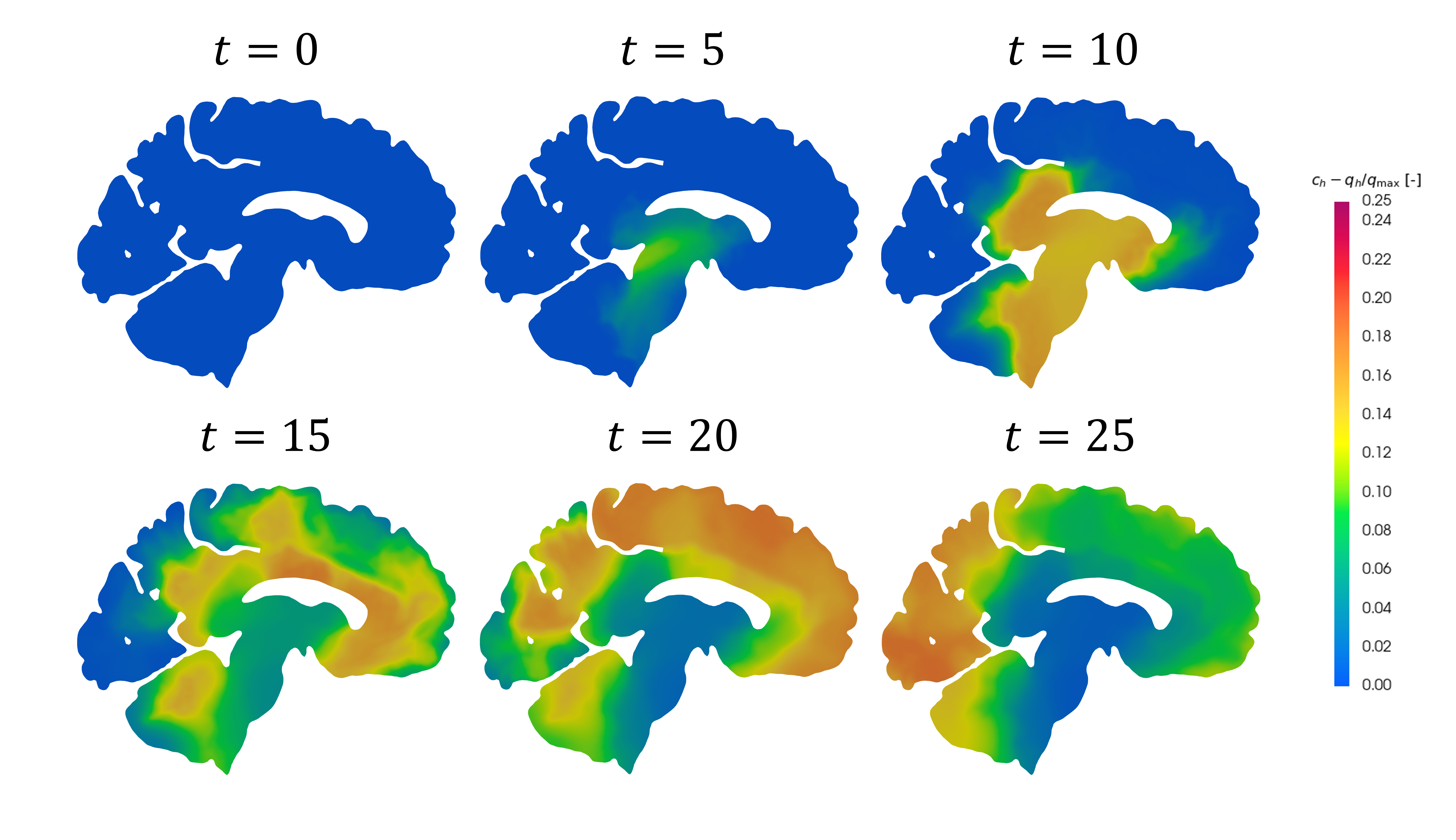}
   \caption{Difference between the solution of the FK model and the solution of the heterodimer one.}
   \label{fig:diff_sol}
\end{figure}

To confirm the fast activation which derives from the application of the FK instead of the heterodimer one, we report the difference of the rescaled solutions in Figure~\ref{fig:diff_sol}. As we can observe, due to the asymptotic behaviour of the two models, the differences are small inside the regions where we do not have misfolded protein concentration $(q_h\simeq 0)$ and inside the ones completely activated $(q_h\simeq q_\mathrm{max})$. On the contrary, where we have a propagating front we can notice a discrepancy in the magnitude with a maximum value of 0.2, which confirms faster dynamics using the FK model. This confirms the differences between the models: the heterodimer model is more informative and an appropriate pathology-specific calibration is needed before evaluating the possibility of reducing it and simulating the FK one.
\par
Finally, concerning the computational costs, we remark that the use of the FK model allows us to solve a system with only $14\,452$ degrees of freedom, half of the $29\,904$ ones required by the heterodimer model. This is reflected in smaller times both in the assembling and solving parts.

\section{Conclusions}
\label{sec:conclusion}
In this work, we discussed the polyhedral discontinuous Galerkin method for the modelling of prionic spreading in brain neurodegeneration, considering in particular a comparison between the heterodimer model and the Fisher-Kolmogorov model.
\par
To carry out the comparison in a realistic configuration we present a simulation of the spreading of tau protein associated with Alzheimer's disease. In particular, we considered a slice of a real brain in the sagittal plane with a polygonal agglomerated grid and on a 3D brain geometry, considering distinct physical parameters in white and grey matters, to describe the domain heterogeneous properties. To compare approximation capabilities with the different models, we compare also the activation times from the simulation results. The obtained results are coherent with the medical literature. However, the results we obtained confirm the need for the heterodimer model whenever the healthy proteins' concentration is far to be constant. We underline that these properties must be analysed specifically for each pathology and prionic protein.

\section*{Acknowledgments}
The brain MRI images were provided by OASIS-3: Longitudinal Multimodal Neuroimaging: Principal Investigators: T. Benzinger, D. Marcus, J. Morris; NIH P30 AG066444, P50 AG00561, P30 NS09857781, P01 AG026276, P01 AG003991, R01 AG043434, UL1 TR000448, R01 EB009352. AV-45 doses were provided by Avid Radiopharmaceuticals, a wholly-owned subsidiary of Eli Lilly.

\section*{Declaration of competing interests}
The authors declare that they have no known competing financial interests or personal relationships that could have appeared to influence the work reported in this article.

\bibliographystyle{hieeetr}
\bibliography{sample.bib}
\end{document}